\documentclass[10pt]{amsart}
  
\usepackage{amscd,graphicx,epsfig}
\usepackage[colorlinks=true, pdfstartview=FitV, linkcolor=blue, citecolor=blue, urlcolor=blue]{hyperref}
\usepackage{amssymb}
\usepackage[all]{xy}
\usepackage{stmaryrd}
\usepackage{mathrsfs}

\newtheorem{thm}{Theorem}
\newtheorem*{thm*}{Theorem}

\newtheorem*{mthm}{Main Theorem}
\newtheorem*{app}{Application}

\newtheorem*{conj*}{Conjecture}

\newtheorem*{prob*}{Problem}
\newtheorem*{satz*}{Satz}

\newtheorem{prop}{Proposition}
\newtheorem*{prop*}{Proposition}
\newtheorem{lem}{Lemma}
\newtheorem*{lem*}{Lemma}

\newtheorem*{cor*}{Corollary}

\theoremstyle{definition}
\newtheorem{defn}{Definition}

\theoremstyle{remark}
\newtheorem*{rem*}{Remark}
\newtheorem{rem}{Remark}

\newcommand{\name}[1]{\textsc{#1\/}}
\newcommand{\NN}{{\mathbb N}}

\newcommand{\CC}{{\mathbb C}}

\renewcommand{\AA}{{\mathbb A}}

\newcommand{\GGG}{\mathcal G}
\newcommand{\UUU}{\mathcal U}

\newcommand{\HHH}{\mathcal H}
\newcommand{\FFF}{\mathcal F}

\newcommand{\VV}{\mathscr V}

\newcommand{\TG}{T{\mathcal G}}

\newcommand{\be}{\begin{enumerate}}
\newcommand{\ee}{\end{enumerate}}

\DeclareMathOperator{\id}{id}

\DeclareMathOperator{\Aut}{Aut}
\DeclareMathOperator{\Lie}{Lie}

\DeclareMathOperator{\Cent}{Cent}

\newcommand{\bbmat}{\begin{bmatrix}}
\newcommand{\ebmat}{\end{bmatrix}}
\newcommand{\bsmat}{\begin{smallmatrix}}
\newcommand{\esmat}{\end{smallmatrix}}

\newcommand{\g}{\mathbf g}

\newcommand{\bu}{\mathbf u}

\newcommand{\f}{\mathbf f}
\newcommand{\e}{\mathbf e}

\newcommand{\bt}{\mathbf t}

\newcommand{\lab}[1]{\label{#1}}

\title[A note on Automorphisms of the Affine Cremona Group]
{A note on Automorphisms of the Affine Cremona Group}
\author{Immanuel Stampfli}

\address{Mathematisches Institut \\
Universit\"at Basel \\ Rheinsprung 21, CH-4051 Basel}
\email{immanuel.e.stampfli@gmail.com}
\thanks{The author is supported by the 
Swiss National Science Foundation (Schweizerischer National\-fonds),
grant ``Automorphisms of Affine $n$-Space" 137679}

\begin{document}

\begin{abstract}
	Let $\GGG$ be an ind-group and let $\UUU \subseteq \GGG$ be a 
	unipotent ind-subgroup.
	We prove that an abstract group automorphism $\theta \colon \GGG \to \GGG$
	maps $\UUU$ isomorphically onto a unipotent
	ind-subgroup of $\GGG$, provided that $\theta$ fixes a closed torus 
	$T \subseteq \GGG$, which normalizes $\UUU$ and the action of 
	$T$ on $\UUU$ 
	by conjugation fixes only the neutral element. As an application
	we generalize a result by \name{Hanspeter Kraft} and the author 
	as follows:
	If an abstract group automorphism of the affine Cremona group $\GGG_3$ in
	dimension 3 fixes the subgroup of tame automorphisms $\TG_3$, 
	then it also fixes a whole family of non-tame automorphisms
	(including the \name{Nagata} automorphism).
\end{abstract}

\maketitle


\setcounter{subsection}{-1}

\subsection{Introduction}

Throughout this note we denote by $\GGG_n$ the group of polynomial
automorphisms $\Aut(\AA^n)$ of the complex affine space $\AA^n = \CC^n$.
Such an automorphism has the form $\g = (g_1, \ldots, g_n) \in \GGG_n$
with polynomials $g_1, \ldots, g_n \in \CC[x_1, \ldots, x_n]$.
We define $\deg \g := \max_i \deg g_i$.
The tame automorphism group $\TG_n$ is the subgroup of $\GGG_n$ 
generated by the affine linear automorphisms (i.e. the automorphisms
$\g$ with $\deg \g \leq 1$) 
and the triangular automorphisms (i.e. the automorphisms $(g_1, \ldots, g_n)$
where $g_i = g_i(x_i, \ldots, x_n)$ depends only on $x_i, \ldots, x_n$ for each $i$).
The main result of \cite{KrSt2013On-Automorphisms-o} is the following.
\begin{thm}	
	\lab{tame.thm}
	Let $\theta \colon \GGG_n \to \GGG_n$ be an abstract automorphism. 
	Then there exist 
	$\g \in \GGG_n$ and a field automorphism $\tau \colon \CC \to \CC$ such that
	\[
		\theta(\f) = \tau(\g \circ \f \circ \g^{-1}) \quad 
		\textrm{for all tame automorphisms $\f \in \TG_n$} \, .
	\]
\end{thm}
If $\theta$ preserves in addition the ind-group structure of $\GGG_n$ 
(see below for a definition), then \name{Alexei Belov-Kanel} and 
\name{Jie-Tai Yu} proved recently that $\theta$
is an inner automorphism of $\GGG_n$ (see \cite{BeYu2013On-The-Zariski-Top}).

In dimension $n=2$ all automorphisms are tame (cf. 
\cite{Ju1942Uber-ganze-biratio} and \cite{Ku1953On-polynomial-ring}).
But in dimension $n=3$, \name{Ivan P. Shestakov} and \name{Ualbai U. Umirbaev} 
showed that the famous \name{Nagata} automorphism $\bu_N \in \GGG_3$ 
(see below for a definition)
is non-tame (cf. \cite{ShUm2004The-tame-and-the-w}). It is an open problem
if there exist non-tame automorphisms in dimension $n > 3$.
A natural question is, whether Theorem~\ref{tame.thm} extends to the entire
automorphism group $\GGG_n$, i.e. whether $\theta(\f) = \tau(\g \circ \f \circ \g)$ 
for all $\f \in \GGG_n$.
If this would be true, then every abstract
automorphism of $\GGG_n$ would preserve the algebraic subgroups of
$\GGG_n$ (see below for a definition). In fact, a main tool in the proof of 
Theorem~\ref{tame.thm} is to show that 
certain algebraic subgroups are sent to isomorphic algebraic subgroups 
under an abstract automorphism of 
$\GGG_n$.
The main point of this note is to refine these techniques. 
In order to state the main result we introduce the concept of an ind-group and
related terms. 

A group $\GGG$ is called an \emph{ind-group} if it is endowed with a filtration by
affine varieties $G_1 \subseteq G_2 \subseteq \ldots \,$, each one closed in the
next, such that $\GGG = \bigcup_{i=1}^{\infty} G_i$ and such that the map
$\GGG \times \GGG \to \GGG$, $(x, y) \mapsto x \cdot y^{-1}$ 
is a morphism of ind-varieties
(see \cite[chap. IV]{Ku2002Kac-Moody-groups-t}  for an introduction to ind-varieties
and ind-groups). We then write $\GGG = \varinjlim G_i$.
For example, $\GGG_n = \varinjlim G_{n, i}$
is an ind-group, where $G_{n, i}$ is the set of all automorphisms 
$\g \in \GGG_n$ with $\deg \g \leq i$ (see \cite{BaCoWr1982The-Jacobian-conje}).
We endow an ind-group $\GGG = \varinjlim G_i$ with the following topology: 
a subset $X \subseteq \GGG$ is closed if and only if $X \cap G_i$ is closed in 
$G_i$ with respect to the Zariski topology for each $i$.
If $\HHH \subseteq \GGG$ is a closed subgroup, then $\HHH$ inherits in a 
canonical way an ind-structure from $\GGG$, namely
$\HHH = \varinjlim \HHH \cap G_i$.
But for our purposes we need a more general definition of an ind-subgroup.

\begin{defn}
	Let $\HHH$ be a subgroup of an ind-group 
	$\GGG = \varinjlim G_i$. 	
	We say that $\HHH$ is an \emph{ind-subgroup} of $\GGG$
	if $\HHH$ can be turned into an ind-group $\HHH = \varinjlim H_k$ such that 
	to every $k$ there exists $i = i(k)$ such that $H_k \subseteq G_i$ is closed. 
	Clearly, the ind-structure of $\HHH$ is then unique.
	We say that $\HHH$ is an \emph{algebraic subgroup} of $\GGG$, if 
	$\HHH$ is closed in $\GGG$ and contained in some $G_i$.
\end{defn}

\begin{defn}
We say that an ind-group $\UUU$ is \emph{unipotent} 
if $\UUU = \varinjlim U_i$ where $U_i$ is a unipotent algebraic group for all $i$. 
\end{defn}

\begin{rem}
Every element in a unipotent ind-group is unipotent. We don't
know whether an ind-group consisting only of unipotent elements
is always unipotent. If the ind-group is commutative,
then we are able to prove this.
\end{rem}

\begin{mthm}
	\lab{main.thm}
	Let $\theta \colon \GGG \to \GGG$ be an abstract automorphism 
	of an ind-group $\GGG$ that is the identity on
	a closed torus $T \subseteq \GGG$.
	If $\UUU \subseteq \GGG$ is a unipotent 
	ind-subgroup that is normalized by $T$ and if
	the neutral element of $\UUU$ is the only element that
	is fixed under conjugation by $T$,
	then $\theta(\UUU)$ is a unipotent ind-subgroup of $\GGG$ and
	$\theta |_\UUU \colon \UUU \to \theta(\UUU)$ is an 
	isomorphism of ind-groups.
\end{mthm}

Recall that there exists a bijective correspondence between locally nilpotent
derivations of $\CC[x_1, \ldots, x_n]$ and unipotent elements of $\GGG_n$, 
given by $D \mapsto \exp(D)$ where
\[
	\exp(D) = \left( \sum_{i=0}^{\infty} \frac{D^i(x_1)}{i!} , \ldots, 
			\sum_{i=0}^{\infty} \frac{D^i(x_n)}{i!} \right)
\]
(see \cite[sec. 1.5]{Fr2006Algebraic-theory-o}).
If $D$ is a locally nilpotent derivation and $f \in \ker D$ then $f D$ is again
a locally nilpotent derivation and we call $\exp(f D)$ a \emph{modification} 
of $\exp(D)$. For example, the \name{Nagata} automorphism 
$\bu_N$ is a modification of $\bu := \exp D$ where
\[
	D = - 2y \frac{\partial}{\partial x} + z \frac{\partial}{\partial y} \ , \quad
	p = xz + y^2 \in \ker D \quad \textrm{and} \quad 
	\bu_N = \exp(p D) \, .
\]
Recently, \name{Shigeru Kuroda} gave a characterization of the non-tame
modifications of certain unipotent automorphisms 
(see \cite[Theorem 2.3]{Ku2011Wildness-of-polyno}).
This result implies that for all 
$f \in \ker D \setminus \CC[z]$ the modification $\exp(f D)$ of $\bu$ is non-tame.
Clearly, all the modifications of $\bu$ lie in the centralizer $\Cent(\bu)$.
As a consequence of our Main Theorem we get the following result.

\begin{app}
	\lab{app.thm}
	Let $\theta: \GGG_3 \to \GGG_3$ be an abstract automorphism
	that is the identity on the tame automorphisms $\TG_3$. 
	Then $\theta$ fixes $\Cent(\bu)$ where $\bu = \exp(D)$
	and $D = -2y \cdot \partial / \partial x +  z\cdot \partial / \partial y$.
	In particular, $\theta$ fixes the non-tame automorphisms
	$\exp(f D)$ where $f \in \ker D \setminus \CC[z]$ and thus 
	$\theta$ fixes the \name{Nagata} automorphism $\bu_N$.
\end{app}

\begin{rem}
	All the results and proves work over any uncountable 
	algebraically closed field of characteristic zero.
\end{rem}

\subsection{Proof of the Main Theorem}
\lab{main.sec}
Let $V$ be a commutative unipotent algebraic group. Recall that $V$ has a unique
$\CC$-vector space structure such that the product in $V$ corresponds to
addition. Also recall that a map of commutative unipotent algebraic 
groups $V \to V'$ is 
a homomorphism of algebraic groups if and only if it is $\CC$-linear.

We start with a lemma that proves the Main Theorem in the case when
$\UUU \subseteq \GGG$ is an algebraic subgroup isomorphic to $\CC^+$.

\begin{lem}
	\lab{unipot.lem}
	Let $\theta: \GGG \to \GGG$ be an abstract automorphism that
	is the identity on a closed torus $T \subseteq \GGG$ and let
	$U \subseteq \GGG$ be an algebraic subgroup isomorphic to $\CC^+$ which
	is normalized by $T$ with character $\lambda$.
	If $\lambda$ is non-trivial, then $\theta(U) \subseteq \GGG$ 
	is an algebraic subgroup
	isomorphic to $\CC^+$ and $T$ normalizes $\theta(U)$ 
	with the same character $\lambda$. 
	Moreover, $\theta |_{U} \colon U \to \theta(U)$ is
	an isomorphism of algebraic groups.
\end{lem}

\begin{proof}
	Let $U' := \theta(U) \subseteq \GGG$. Choose $\bu_0 \in U$
	that is different from the neutral element $\e \in \GGG$. Then, 
	$U' \setminus \{ \e \} = \{ \, \bt \cdot \theta(\bu_0) \cdot \bt^{-1} \ | \ \bt \in T \, \}$
	and $\{ \e \}$ are constructible subsets of some filter set of $\GGG$ 
	and since $U'$ is a group it follows that $U' \subseteq \GGG$ is an algebraic
	subgroup (see \cite[7.4~Proposition~A]{Hu1975Linear-algebraic-g}). 
	Since $U'$ has no element $\neq \e$ of finite order, $U'$
	is unipotent. As $U'$ is a toric variety with exactly two orbits, 
	$U'$ is one-dimensional 
	(see also \cite[Proposition~2]{KrSt2013On-Automorphisms-o}).
	Let $\lambda'$ be the character of $U'$. We have
	\begin{equation}
		\label{eq}
		\tag{$\ast$}
		\theta(\lambda(\bt) \bu_0) = 
		\theta(\bt \cdot \bu_0 \cdot \bt^{-1}) = 
		\bt \cdot \theta(\bu_0) \cdot \bt^{-1} 
		= \lambda'(\bt) \theta(\bu_0) \quad \textrm{for all} \ \bt \in T \, .
	\end{equation}
	Hence, it follows that $\lambda$ and $\lambda'$ have the same kernel
	and thus $\lambda = \pm \lambda'$. If we take $\bt \in T$
	such that $\lambda(\bt) = 2$ then eq.~\eqref{eq} implies that
	$\lambda' \neq -\lambda$. Hence, $\lambda = \lambda'$ and
	$\theta |_{U} \colon U \to U'$ is $\CC$-linear by eq.~\eqref{eq}.
\end{proof}

\begin{proof}[Proof of the Main Theorem]
	Let $U \subseteq \UUU$ be an algebraic subgroup that is normalized by $T$.
	Choose closed algebraic subgroups $V_1, \ldots, V_r \subseteq U$
	which are isomorphic to $\CC^+$ and which are normalized by $T$, 
	such that $\Lie U = \Lie V_1 \oplus \ldots \oplus \Lie V_r$.
	Thus, for suitable indices $i_1, \ldots, i_m$ we
	have $U = V_{i_1} \cdot \ldots \cdot V_{i_m}$ 
	(see \cite[7.5~Proposition]{Hu1975Linear-algebraic-g}). 
	According to Lemma~\ref{unipot.lem}, 
	$\theta(U) = \theta(V_{i_1}) \cdot \ldots \cdot \theta(V_{i_m})$ is a constructible
	subset of some $G_i \subseteq \GGG$ and thus $\theta(U)$ is an algebraic 
	subgroup of $\GGG$. Again, since no element $\neq \e$ in $U$ has
	finite order, $\theta(U)$ is unipotent.
	Consider the following commutative diagram.
	\[
		\xymatrix{
			V_{i_1} \ar@{->>}[d] \times \ldots \times V_{i_m} 
			\ar[rr]^-{\theta \times \ldots \times \theta} && 
			\theta(V_{i_1}) \times \ldots \times \theta(V_{i_m}) \ar@{->>}[d] \\
			U \ar[rr]^-{\theta |_U} && \theta(U)
		}
	\]
	The vertical maps are induced by the product in $\UUU$ and hence 
	they are surjective
	morphisms. The top horizontal map is 
	an isomorphism of varieties by Lemma~\ref{unipot.lem}. 
	The lemma below due to 
	\name{Hanspeter Kraft} implies that the abstract group homomorphism 
	$\theta |_U$ is an isomorphism of algebraic groups.

	As we can replace 
	the filtration of $\UUU$ by a filtration of unipotent algebraic subgroups,
	each one normalized by $T$, it follows that 
	$\theta(\UUU) \subseteq \GGG$
	is a unipotent ind-subgroup and 
	$\theta |_\UUU \colon \UUU \to \theta(\UUU)$ is an 
	isomorphism of ind-groups.
\end{proof}

\begin{lem}[\name{Hanspeter Kraft}]
	Let $X$ and $Y$ be affine varieties and let 
	$f \colon X \to Y$ be an abstract map. If there exists a surjective morphism
	$g \colon Z \twoheadrightarrow X$ such that the composition 
	$f \circ g \colon Z \to Y$ is a morphism and if $X$ is normal, 
	then $f$ is a morphism.
\end{lem}

For the proof of this lemma, one shows that the graph of $f$
is closed in $X \times Y$.

\subsection{Proof of the Application}
\label{app.sec}
First, we determine the structure of $\Cent(\bu)$. 
Denote by $E$ the partial derivative 
with respect to $x$. The ind-subgroups of $\GGG_3$
listed below are clearly contained in $\Cent(\bu)$.
\begin{eqnarray*}
	C & := & \{ \, (a x, a y, a z) \ | \ a \in \CC^\ast \, \} \\
	\FFF & := & \{ \, \exp(fD) \ | \ f \in \ker D \, \} \\
	\HHH & := & \{ \, \exp(hE) \ | \ h \in \ker E \cap \ker D \, \}
\end{eqnarray*}
Remark that $\ker E \cap \ker D = \CC[z]$ and thus
$\HHH = \{ \, (x + h, y, z) \ | \ h \in \CC[z] \,  \}$.

\begin{prop}
	\lab{cent.prop}
	We have a semi-direct product decomposition 
	\[
		\Cent(\bu) = C \ltimes (\HHH \ltimes \FFF) \, .
	\]
\end{prop}

\begin{proof}
	Recall that $p = xz + y^2$, $z \in \ker D$. In fact, 
	$R := \ker D$ is the polynomial ring $\CC[z, p]$ by 
	\cite[Theorem~(b)]{Ba1984A-nontriangular-ac}. We have 
	$R[x] = \CC[z, x, y^2]$ and hence a decomposition
	$\CC[x, y, z] = R[x] \oplus y R[x]$.
	Let $\g = (g_1, g_2, g_3) \in \Cent(\bu)$.
	Write $g_1 = v + y q$ with polynomials $v, q \in R[x]$.
	In $\CC[x, y, z, t]$ we have, by definition,
	\[
		v(x - 2ty - t^2 z) + (y + tz)q(x - 2ty -t^2 z) \ = \ 
		v(x) + yq(x) - 2t g_2 -t^2 g_3 \, .
	\]
	A comparison of the coefficients with respect to the variable $t$ shows that
	$v = r + sx$ with $r, s \in R$, and $q \in R$. Hence, we get
	$g_1 = r + sx + qy$, $g_2 = sy - (q/2)z$, $g_3 = sz$ and $s \in \CC^\ast$. 
	Up to post composition with an element of $C$
	we can assume that $s = 1$. Thus,
	\[
		\g \circ \exp \left( \frac{q}{2} D \right)
		= (x + r + \frac{q^2}{4} z, y, z) \, .
	\]
	One easily sees that this automorphism belongs to $\HHH$.
\end{proof}

Let  $\UUU$ be the ind-subgroup $\HHH \circ \FFF \subseteq  \GGG_3$. 
Every element $\g = \exp(hE) \circ \exp(fD) \in \HHH \circ \FFF$ satisfies
\[
	\g^m = \exp(mhE) \circ \exp\left( \sum_{i = 0}^{m-1} \exp(ihE)^{\ast}(f) D \right) 
	\quad \textrm{for all $m \in \NN$} \, .
\]
Since $\HHH \cap \FFF = \{ \bf{id} \}$, 
this shows that every element $\neq \bf{id}$ of $\UUU$ has infinite order and is contained
in an algebraic subgroup of $\UUU$. Thus every element of $\UUU$ is unipotent.
Now, Proposition~\ref{cent.prop} implies that $\UUU$ is the set of unipotent 
elements of $\Cent(\bu)$. A calculation shows that
\[
	\label{triangle}
	\tag{$\triangle$}
	\FFF = \Cent_{\UUU} [\UUU, \UUU]
\]
where $[\UUU, \UUU]$ is the commutator subgroup of $\UUU$.
Let $T := \{ \, (a^2 b^{-1} x, a y, b z)  \ | \ a, b \in \CC^\ast \, \}$ which is a closed 
algebraic subgroup of $\GGG_3$. The torus $T$ normalizes 
$\Cent(\bu)$, $\HHH$, $\FFF$ and $\UUU$.
In fact, it follows from eq. \eqref{triangle} that $T$ is the 
largest subgroup of the standard torus in $\GGG_3$ that normalizes 
$\Cent(\bu)$.

\begin{proof}[Proof of the Application]
	As $\bu$ is a tame automorphism, $\theta$ preserves $\Cent(\bu)$. 
	A calculation shows that the neutral element
	is the only element of $\HHH$ (of $\FFF$) that commutes with $T$.
	Thus the Main Theorem applied to 
	the unipotent ind-subgroups $\HHH$ and $\FFF$ of $\GGG_3$ 
	implies that $\theta(\UUU) = \UUU$. 
	From eq. \eqref{triangle} it follows that $\theta$ preserves 
	$\FFF$. We define $\VV_{\FFF}$ as the set of all algebraic 
	subgroups of $\FFF$ 
	which are isomorphic to $\CC^+$ and which are normalized by $T$.
	One can see that the elements of $\VV_\FFF$ 
	correspond to the locally nilpotent derivations 
	$z^m p^n D$ for all $m, n \geq 0$ and hence different elements in
	$\VV_\FFF$ have different characters. This implies that for all 
	$V \in \VV_{\FFF}$ there exists $a_V \in \CC^\ast$, such that
	$\theta |_V = a_V \id_V$. Let $\e := (x-1, y, z) \in \HHH$. A 
	calculation shows that
	\begin{equation*}
		\label{eq2}
		\e \circ \exp(z^m p^n D) \circ \e^{-1} = 
		\exp \left(\sum_{i = 0}^n {n \choose i} z^{i+m} p^{n-i} D \right)  
		\quad \textrm{for all $m, n \in \NN$}\, .
	\end{equation*}
	Applying $\theta$ to the last equation and using the fact that 
	$\exp(z^j D)$ and $\e$ are tame automorphisms yields
	$a_V = 1$ for all $V \in \VV_\FFF$. As all the subgroups 
	$V \in \VV_{\FFF}$ generate $\FFF$, it follows that $\theta$ is the identity
	on $\FFF$. Since $\HHH$ and $C$ consist of tame automorphisms
	this finishes the proof.
\end{proof}

\section*{Acknowledgements}
I would like to thank \name{Hanspeter Kraft} for many fruitful discussions.

\vskip 0.2cm

\bibliographystyle{amsalpha}
\bibliography{NoteAutAnBIB}

\vskip.2cm

\end{document}